# Examining the Impact of Tutorial Activity Engagement on Undergraduate Students' Collaborative Preferences

Sang Hyun Kim
*University of Auckland*
skim660@aucklanduni.ac.nz

Tanya Evans
*University of Auckland*
t.evans@auckland.ac.nz

This study examines the impact of tutorial engagement on Collaborative Preferences for Learning Mathematics (CPLM) in a tertiary context. A two-way mixed ANOVA analysed these preferences over a semester in a sample of undergraduate students. As expected, collaborative engagement had a significant main effect, with students who collaborated more reporting stronger preferences for working with their peers (higher CPLM). The absence of an interaction effect between the nature of tutorial engagement and time suggests CPLM differences remain stable. This may indicate that familiar modes of tutorial engagement may reinforce existing collaboration preferences.

Collaboration is a central component of mathematics learning, particularly in tertiary education, where collaborative and active learning approaches are increasingly emphasised, and students differ in their preferences for working with others. These preferences, termed Collaborative Preferences for Learning Mathematics (CPLM), reflect the underlying cognitive, social, and motivational needs and dispositions that shape how students engage in mathematical tasks (Kim & Evans, in press). While extensive research has explored strategies that promote collaboration (e.g., Leitner & Gabel, 2024; Lo et al., 2017; Oates et al., 2024; Oates et al., 2016), less is known about how students' engagement with mathematics influences their preferences for working with peers over time. Although existing studies highlight both the benefits and challenges of collaborative learning, the role of tutorial participation in shaping these preferences remains underexplored.

This study examines the role of tutorial engagement—whether students engage primarily individually, in a mixed format, or collaboratively—in shaping CPLM. Using a two-way mixed ANOVA, we investigate changes in CPLM across a semester, exploring whether tutorial participation fosters shift in students' collaborative learning preferences. Specifically, we address two research questions: (1) What are the CPLMs among undergraduate students based on tutorial engagement, and (2) how does tutorial engagement (individual learners, mixed learners, collaborative learners) influence students' CPLM over time? By understanding these dynamics, this research contributes to ongoing discussions in mathematics education regarding the interplay between instructional contexts and students' learning preferences.

## Literature Review

### Collaborative Preferences

Preferences students have within the domain of mathematics are deeply connected to the values and interests they uphold. These inclinations to work in particular ways may reflect underlying cognitive and affective needs that stem from their past experiences and motivate behaviour and effort expended by the students (Goldin et al., 2016; Saadati & Reyes, 2019). As a driver of behaviour of short- and long-term habitual practices, student preferences are likely to have an influence on the ways students engage in mathematical thinking and communication and, hence, have an indirect impact on their learning outcomes. One key dimension of these preferences is their disposition toward collaborative versus individual learning environments. While some students thrive in collaborative settings, valuing the opportunity to share diverse perspectives, others prefer solitary study to maintain focus and control over their learning process.

Mathematics education research provides further insight into the complexities of students' collaborative behaviours and learning preferences. Studies examining undergraduate mathematics students have shown that, while many hold positive attitudes toward group work, they often still prefer independent learning (MacBean et al., 2004). Informal group interactions outside of the classroom offer flexibility in how students can collaborate, and the unsupervised nature of it can foster a sense of camaraderie as they share experiences in problem-solving (MacBean et al., 2004). Formal group activities taking place during class hours can be valuable as instructors can have greater oversight over the students, allowing them to monitor student progress and offer guidance. Nevertheless, they also present challenges, including unequal participation and difficulties working with unfamiliar peers (D'Souza & Wood, 2003; Karau & Williams, 1993; MacBean et al., 2004; North et al., 2000). Further illustrating these complexities, a study conducted in an advanced-level mathematics course found that students perceived the greatest benefit of group work to be the exposure to diverse perspectives (Sheryn & Ell, 2014). Not all students viewed group work positively. Many who were initially apprehensive about it later described the experience as "interesting," even if they did not necessarily find it useful (Sheryn & Ell, 2014, p. 875). D'Souza and Wood (2003) reported on the cognitive advantages of learning with peers, such as fostering higher-level thinking skills and improving the performance of weaker students when paired with higher-achieving peers. These benefits, however, were accompanied by challenges, such as disengaged participants, which resulted in more engaged peers expressing frustration. Research in broader educational contexts echoes these concerns, noting that while collaboration can enhance learning, unproductive group dynamics and inefficiencies are recurring challenges (Salomon & Globerson, 1989; Sofroniou & Poutos, 2016; Tucker & Abbasi, 2016).

Considering these insights, the decision to incorporate group work requires careful thought about both student preferences and class dynamics. This complexity highlights that students' preferred learning environments do not always align with those that best support their learning outcomes. For example, Deslauriers et al. (2019) found that students often misjudge the efficacy of different learning approaches, favouring methods that create an illusion of learning, which may not always lead to greater understanding. Additionally, while early research suggests that aligning instructional practices with student preferences may enhance engagement and learning in certain contexts, this approach requires careful consideration (Okebukola, 1986).

Beyond immediate educational contexts, research in psychology and education has identified patterns between students' collaborative preferences and individual traits, including gender and personality, suggesting that these preferences may be shaped to some extent by innate factors. For instance, Owens (1985) found that girls generally exhibit a stronger preference for cooperative learning than boys, with this tendency becoming more pronounced as they age. Regarding personality traits, Chamorro-Premuzic et al. (2007) observed that extroverted students are more likely to favour collaborative learning, while those with higher neuroticism scores tend to prefer independent study. While these traits may influence students' predispositions, they do not define them, as students' preferences for collaboration can also be shaped by their broader social, cognitive, and educational experiences.

**Theoretical Grounding**

The familiarity principle, also known as the mere-exposure effect, posits that individuals tend to develop a preference for stimuli they encounter frequently (Zajonc, 1968; Zajonc, 2001). In educational contexts, this suggests that repeated exposure to specific instructional practices can influence students' preferences and perceptions of their effectiveness. For instance, when students consistently experience collaborative learning environments, they may develop a stronger preference for such settings due to increased familiarity and comfort (Macaluso et al., 2022). This repeated exposure can lead to increased comfort and perceived efficacy in group settings, thereby reinforcing positive attitudes toward collaborative learning. Conversely, students who predominantly experience individualised learning environments may become more comfortable with solitary study methods, thereby reinforcing preferences for independent work. This alignment between exposure and preference reflects the broader influence of familiarity on perceived learning effectiveness. A study by Vroom et al. (2022) provides additional evidence of the mere-exposure effect within the context of mathematics education. Their findings show that students often perceive instructional practices they encounter frequently as more beneficial to their learning. Specifically, students in this study reported that interactive classroom characteristics, such as group work and peer support, were helpful when experienced on a regular basis. These findings emphasise the role of familiarity in shaping students' attitudes towards particular teaching methods and suggest that the repeated use of collaborative or individual approaches can influence how students evaluate their effectiveness.

In the context of this study, the familiarity principle is particularly relevant for examining how students' preferences for collaborative or individual learning in mathematics might develop based on their cumulative instructional experiences. By exploring the relationship between repeated exposure to specific learning modes and students' collaborative learning preferences, this research aims to contribute to a deeper understanding of how instructional practices can shape affective factors, such as engagement, comfort, and motivation, in mathematics education. This theoretical framing provides a foundation for the research questions and aims, which seek to investigate how exposure to collaborative or individual learning influences students' preferences and whether these preferences relate to academic performance and broader affective outcomes.

## Data

The data used in this study were collected during the second semester of 2024 from a second-year undergraduate mathematics course at the University of Auckland. The course covered three primary topics: Calculus II, Linear Algebra II, and Differential Equations. Over the 12-week semester, students were expected to engage in weekly one-hour tutorials focused on problem-solving with a tutor —either the lecturer or an experienced graduate student—who was available to help students with their progress. Students were encouraged, though not required, to collaborate with peers. This study employed a quantitative, survey-based methodology to examine students' collaborative learning preferences and their self-reported tutorial engagement. Data collection involved two primary measures: (1) a Collaborative Preference for Learning Mathematics (CPLM) scale, which assessed students' preferences for individual versus collaborative learning, and (2) a self-reported measure of tutorial engagement, which captured the extent to which students worked individually or collaboratively throughout the semester. The self-reported CPLM data and tutorial engagement data were collected at three points in the semester (the first week, after the midsemester break, and the final week). Of the 294 students enrolled in the course, 201 completed the surveys at all three time points, and these responses were used in the analysis.

A validated 5-item scale was used to capture students' CPLM (Kim & Evans, in press). This scale captures students' preferences for collaborative versus individual learning across a range of scenarios. Each item was presented as a slider-based, close-ended question, with responses recorded on a scale from 0 (indicating a preference for individual learning) to 100 (indicating a preference for collaborative learning). The five items were designed to evaluate preferences across different contexts, including learning new concepts, studying for assessments, and engaging in problem-solving activities. The prompt for the scale was: "Consider yourself learning mathematics. By moving the slider, state the extent to which you prefer to do it from 0 = Individually to 100 = Collaboratively for each of the items." The five items were as follows: (1) "What is the most effective way for you to learn mathematics?" (CPLM_1), (2) "What is the best way for you to make sense of mathematics?" (CPLM_2), (3) "What is the most effective way for you to study for high-stakes maths assessments (e.g., exams)?" (CPLM_3), (4) "In what social setting do you prefer to be exposed to novel concepts?" (CPLM_4), and (5) "In what social setting do you prefer to engage in problem-solving in low-stakes assessment (e.g., homework, practice exercises)?" (CPLM_5). To quantify CPLM, a composite score was calculated for each student at each time point by averaging their responses across all five items. This approach ensured that students' collaborative learning preferences were captured holistically rather than being based on a single item. Averaging across multiple items reduced variability and improved reliability, providing a more stable measure of individual differences in learning preferences.

Collaborative tutorial engagement was assessed through a self-reported survey item designed to capture students' actual behaviours during tutorials over the semester. Students responded to the following prompt: "Please indicate below to what extent you worked individually or with others in your tutorials during the entire semester." Responses were recorded on a five-point Likert scale: 1 = Always individually ($n = 18$), 2 = Mostly individually ($n = 38$), 3 = Sometimes individually and sometimes with others ($n = 61$), 4 = Mostly with others ($n = 41$), and 5 = Always with others ($n = 43$). To ensure more comparable group sizes for analysis, responses were recoded into three broader categories: individual learners, comprising students who responded "Always individually" or "Mostly individually" ($n = 56$); mixed learners, comprising students who responded "Sometimes individually and sometimes with others" ($n = 61$); and collaborative learners, comprising students who responded "Mostly with others" or "Always with others" ($n = 84$). This recoding facilitated meaningful statistical comparisons while preserving the integrity of students' self-reported tutorial engagement.

## Analysis

Descriptive statistics were reported for students' CPLM scores, and a two-way mixed ANOVA was conducted to examine how CPLM scores changed over time and whether this change varied based on students' tutorial engagement. The between-subjects factor was collaborative tutorial engagement, with three levels (individual, mixed, and collaborative learners), while the within-subjects factor was time, measured at three points during the semester. An interaction term was included to assess whether the change in CPLM scores over time differed depending on tutorial engagement style.

A two-way mixed ANOVA was appropriate for this analysis as CPLM scores were measured on a continuous scale, allowing for the comparison of mean differences both within and between groups while accounting for repeated measures. As the test assumes sphericity for repeated measures, Mauchly's test of sphericity was conducted to verify this assumption. The test revealed a violation to the sphericity assumption for the two-way interaction, $\chi^2(2) = 9.654$, $p = 0.008$. To address this, Greenhouse-Geisser corrections ($\varepsilon = .954$) were applied to adjust the degrees of freedom for the F-tests.

## Results

### Descriptive Statistics

Table 1 presents the descriptive statistics for CPLM scores based on the recoded tutorial engagement groups—individual learners, mixed learners, and collaborative learners. Individual learners, who primarily worked alone in tutorials, consistently reported the lowest mean CPLM scores across all three time points, with a decline from 43.68 at Time 1 to 38.36 at Time 3. Conversely, collaborative learners, who predominantly worked with others, demonstrated the highest mean scores at each time point, increasing from 53.10 at Time 1 to 56.03 at Time 3, reflecting a strengthening preference for collaborative learning over time. Mixed learners, who engaged in both individual and collaborative approaches, exhibited more dynamic changes, with their mean CPLM scores rising from 46.08 at Time 1 to 51.56 at Time 3. These trends suggest that students who engaged more collaboratively in tutorials increasingly favoured collaborative learning, while those who worked individually showed a decline in collaborative preference over the semester.

**Two-Way Mixed ANOVA**

The interaction between collaborative tutorial engagement and time approached statistical significance but did not meet the conventional threshold, $F(3.817, 377.926) = 2.127$, $p = 0.08$, partial $\eta^2 = 0.02$. This suggests a potential, but not definitive, difference in how CPLM scores changed over time depending on tutorial engagement style.

The main effect of time was not statistically significant, $F(1.909, 377.926) = 1.066$, $p = 0.343$, partial $\eta^2 = 0.005$, indicating that overall CPLM scores did not exhibit a uniform change across the semester. However, the main effect of collaborative tutorial engagement was statistically significant, $F(2, 198) = 8.901$, $p < 0.001$, partial $\eta^2 = 0.082$, indicating a medium effect size (Cohen, 2003). This finding suggests that students' preferences for collaborative learning in mathematics differed significantly depending on their tutorial engagement style, with collaborative learners consistently reporting the strongest preference for collaboration.

**Table 1**

*CPLM Descriptive Statistics at each Time Point Based on Collaborative Tutorial Engagement (three levels)*

| Time | Collaborative tutorial engagement | Mean | Std. deviation | N |
|---|---|---|---|---|
| 1 | Individual Learners | 43.68 | 22.48 | 56 |
|   | Mixed Learners | 46.08 | 20.69 | 61 |
|   | Collaborative Learners | 53.10 | 22.95 | 84 |
|   | Total | 48.65 | 22.41 | 201 |
| 2 | Individual Learners | 39.71 | 24.90 | 56 |
|   | Mixed Learners | 46.47 | 19.90 | 61 |
|   | Collaborative Learners | 53.56 | 21.62 | 84 |
|   | Total | 47.55 | 22.72 | 201 |
| 3 | Individual Learners | 38.36 | 25.19 | 56 |
|   | Mixed Learners | 51.56 | 19.60 | 61 |
|   | Collaborative Learners | 56.03 | 21.31 | 84 |
|   | Total | 49.75 | 23.07 | 201 |

**Figure 1**

*Change in CPLM Scores Over Time for Individual, Mixed, and Collaborative Learners*

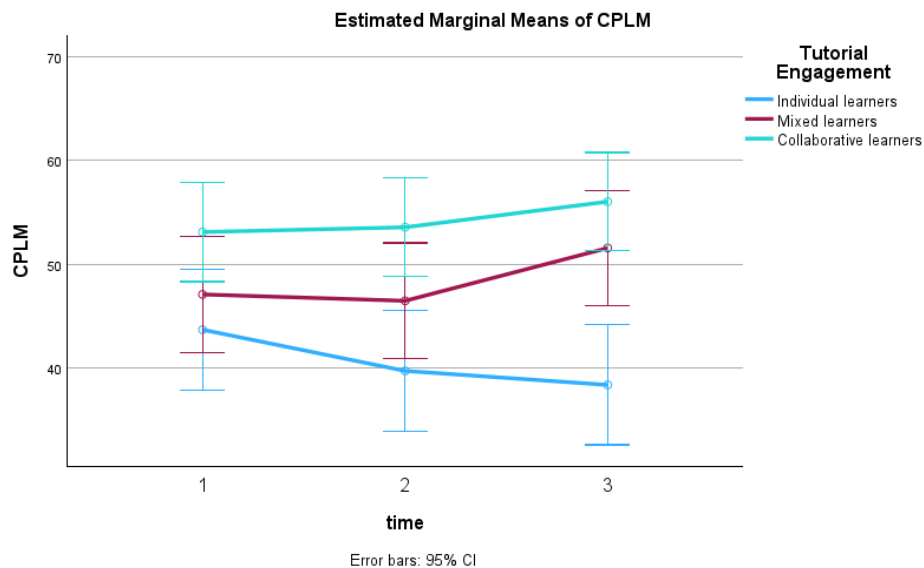

## Discussion

This study investigated how collaborative tutorial engagement influences students' collaborative preferences for learning mathematics (CPLM) over time. The findings reveal that students who engaged collaboratively reported higher CPLM scores and exhibited a positive upward trajectory, whereas individual learners experienced a steady decline in their collaborative preferences. These results suggest that educational contexts and structured peer interactions may play a pivotal role in shaping and reinforcing students' learning preferences over time. Regarding the second research question, while no significant interaction between collaborative tutorial engagement and time emerged, the main effect of collaborative tutorial engagement was significant. Specifically, students who actively worked with peers throughout the semester demonstrated higher preferences for collaboration by the end of the semester compared to those who predominantly engaged in individual learning.

One interpretation of this result is that the nature of students' interactions and engagement during tutorials may reinforce beliefs and preferences that align with their experiences. This interpretation aligns with the familiarity principle (Zajonc, 1968; Zajonc, 2001), which suggests that individuals develop preferences for contexts and experiences they encounter frequently. Regular exposure to peer interactions and collaborative problem-solving may have fostered a sense of comfort, confidence, and perceived efficacy in collaborative settings for students who participated in group activities. Conversely, individual learners, with less exposure to collaborative opportunities, may have increasingly relied on solitary learning strategies, thereby reinforcing their preferences for independent study. These dynamics highlight the importance of considering how repeated exposure to specific educational practices can shape students' learning preferences and behaviours over time.

This finding is consistent with prior research, such as Vroom et al. (2022), who found that students often perceive the instructional practices they encounter most frequently as the most beneficial for their learning. Specifically, Vroom et al. (2022) noted that students viewed classroom characteristics as helpful when they facilitated deeper conceptual understanding, built on prior knowledge, supported its application in new contexts, and provided opportunities to explain concepts to peers. In the present study, students engaged in collaborative tutorials may have experienced similar benefits through interactions with peers, which could have reinforced their preference for collaborative learning. On the other hand, students who predominantly engaged in individual learning may have found value in refining independent strategies that aligned with their preferred learning style, strengthening their inclination towards solitary study.

Broader implications of collaborative learning preferences should also be considered. Students who prefer collaborative learning may benefit from shared problem-solving experiences and opportunities to enhance social and communication skills. Conversely, students who prefer individual study may develop strengths in self-regulation, independent thinking, and autonomous learning. This study does not advocate for one preference over another, nor does it propose definitive recommendations for practice. Instead, it emphasises the importance of understanding how learning preferences are shaped by repeated behaviours and exposure to particular instructional contexts. Future research should explore the long-term implications of these preferences, particularly their impact on academic outcomes, educational trajectories, and students' capacity for lifelong learning. Understanding how preferences evolve in different learning environments and how they influence engagement, achievement, and broader skill development may further enhance educational practices aimed at supporting diverse learners.

## Acknowledgements

Ethics approval UAHPEC26081 was granted by the University of Auckland Human Participants Ethics Committee on 5/07/2023 for three years, and participants gave informed consent for their data to be used for research purposes.